\documentstyle{amsppt}
\hoffset=0.8in
\pagewidth{125mm}
\pageheight{195mm}
\parindent=8mm
\frenchspacing \tenpoint
\TagsAsMath
\NoRunningHeads
\topmatter
\title
Weighted Strichartz estimates and global existence for
semilinear wave equations
\endtitle
\author 
Vladimir Georgiev, Hans Lindblad and
Christopher D. Sogge
\endauthor
\thanks  The first author was partially supported by the 
Alexander von Humboldt Foundation and contract MM-516 with
the Bulgarian Ministry of Education, Science and Technology.
The last two authors were supported in part by the National 
Science Foundation.  \endthanks
\address
Bulgarian Academy of Sciences 
\endaddress
\address
University of California, San Diego
\endaddress
\address University of California, Los Angeles
\endaddress
\endtopmatter

\define\st{\Bbb R^{1+n}_+}
\define\Rn{\Bbb R^n}

\redefine\epsilon{\varepsilon}

\define\cd{\, \cdot \, }

\define\dist{\text{{\rm dist }}}

\define\supp{\text{{\rm supp }}}

\widestnumber\key{999}

\define\wght{(t^2-|x|^2)}
\define\whgt{(t^2-|x|^2)}

\document

\head {\bf 1. Main results} \endhead

The purpose of this paper is to prove sharp global existence theorems
in all dimensions for small-amplitude wave equations with power-type
nonlinearities.  For a given ``power'' $p>1$, we shall therefore consider
nonlinear terms $F_p$ satisfying
$$\bigl| \, (\partial/\partial u)^jF_p(u)\, \bigr|\le C_j|u|^{p-j}, \, \,
j=0,1.
\tag1.1$$
The model case, of course, is $F_p(u)=|u|^p$.  If $\st=\Bbb R_+
\times \Rn$, and if $f,g\in C^\infty_0(\Rn)$ are fixed, we shall consider
Cauchy problems of the form
$$\cases
\square u=F_p(u), \, \, \, (t,x)\in \st
\\
u(0,x)=\epsilon f(x), \, \, \partial_tu(0,x)=\epsilon g(x),
\endcases
\tag1.2$$
where $\square =\partial^2/\partial t^2-\Delta_x$ denotes the D'Alembertian.
Our chief goal then is to find, for a given $n$, the sharp range of powers
for which one always has a global weak solution of (1.2) if $\epsilon>0$ is
small enough.

Note that, even in the linear case, where one solves an inhomogeneous 
equation with a Lipschitz forcing term, in general one can only obtain weak
solutions.  An interesting problem would be to find out to what degree the
regularity assumptions on the data can be relaxed in the spirit of
\cite{8}; however, we shall not go into that here.

Let us now give a bit of historical background.
In 1979, John \cite{6} showed that when $n=3$ global solutions always exist
if $p>1+\sqrt2$ and $\epsilon>0$ is small.  He also showed that the power
$1+\sqrt2$ is critical in the sense that no such result can hold if $p<1+\sqrt2$
and $F_p(u)=|u|^p$.  It was shown sometime later by Schaeffer \cite{12} that
there can also be blowup for arbitrarily small data in $(1+3)$-dimensions when
$p=1+\sqrt2$.

The number $1+\sqrt2$ appears to have first arisen in Strauss' work \cite{21}
on scattering for small-amplitude semilinear Schr\"odinger equations.  Based
on this, he made the insightful conjecture in \cite{22} that when $n\ge2$ 
global solutions of (1.2) should always exist if $\epsilon$ is small and $p$
is greater than a critical power which is the solution of the quadratic equation
$$(n-1)p_c^2-(n+1)p_c-2=0, \, \, \, p_c>1.
\tag1.3$$

This conjecture was shortly verified when $n=2$ by Glassey \cite{3}.  John's
blowup results were then extended by Sideris \cite{15}, showing that, for all 
$n$, there can be blowup for arbitrarily small data if $p<p_c$.  In the other direction, 
Zhou \cite{26} showed that when $n=4$, in which case $p_c=2$, there is always
global existence for small data if $p>p_c$.  This result has recently been extended
to dimensions $n\le 8$ in Lindblad and Sogge \cite{9}.  Here it was also shown
that, under the assumption of spherical symmetry, for arbitrary $n\ge3$ global
solutions of (1.2) exist if $p>p_c$ and $\epsilon$ is small enough.  For odd
spatial dimensions, the last result was obtained independently by Kubo \cite{7}.

In this paper we shall show that the assumption of spherical symmetry can be
removed.  Specifically, we have the following

\proclaim{Theorem 1.1}  Let $n\ge3$ and assume that $F_p$ satisfying $(1.1)$
is fixed with $p_c<p\le (n+3)/(n-1)$.  
Then if $\epsilon>0$ is sufficiently small $(1.2)$ has
a unique (weak) global solution $u$ verifying 
$$(\, 1+|t^2-|x|^2|\, )^\gamma u\in L^{p+1}(\st),
\tag1.4$$
for some $\gamma$ satisfying
$$1/p(p+1)<\gamma<(\, (n-1)p-(n+1)\, )/2(p+1).
\tag1.5$$
\endproclaim

Note that our condition on $\gamma$ only makes sense if $p>p_c$.  For, by
(1.3), $1/p(p+1)<((n-1)p-(n+1))/2(p+1)$ if and only if $p>p_c$.

In Theorem 1.1 we have only considered powers smaller than the conformally
invariant power $p_{\text{conf}}=(n+3)/(n-1)$ since it was already known
that there is global existence for powers larger than $p_{\text{conf}}$.
See, e.g., \cite{8}.

We shall prove Theorem 1.1 using certain ``weighted Strichartz estimates''
for the solution of the linear inhomogeneous wave equation
$$\cases
\square w(t,x)=F(t,x), \quad (t,x)\in \st
\\
0=w(0,\cd)=\partial_tw(0,\cd).
\endcases
\tag1.6$$
This idea was initiated by Georgiev \cite{2}.

Before stating our new estimates, though, let us recall the approach that
John \cite{6} used to show that there is global existence for (1.2) when
$n=3$, $p>1+\sqrt2$ and $\epsilon$ is small.  The main step in his proof of
this half of his theorem was to establish certain pointwise estimates for
the solution of (1.6).  Specifically, he proved an inequality which is
equivalent to the following:
$$\multline
\|t(t-|x|)^{p-2}w\|_{L^\infty(\Bbb R^{1+3}_+)}
\le C_p\|t^p(t-|x|)^{p(p-2)}F\|_{L^\infty(\Bbb R^{1+3}_+)}, 
\\
\text{if } \, \, F(t,x)=0, \, \, t-|x|\le1, \, \,
\text{and } \, \, 1+\sqrt2<p\le 3.
\endmultline$$
Since the powers of the weights behave well with respect to iteration, it is
easy to show that this inequality implies that global solutions of (1.2)
exist when $n=3$ if the data is small and $1+\sqrt2<p\le 3$ (cf. Lemma 1.3
below).

Unfortunately, no such pointwise estimate can hold in higher dimensions due to
the fact that fundamental solutions for $\square$ are no longer measures when
$n\ge4$.  Despite this, it turns out that certain estimates, involving simpler
weights which are invariant under Lorentz rotations, hold if one is willing to
consider dual spaces.  Specifically, we have the following

\proclaim{Theorem 1.2}  Suppose that $n\ge2$ and that $w$ solves the linear
inhomogeneous wave equation $(1.6)$ where the forcing term is assumed to
satisfy $F(t,x)=0$ if $t-|x|\le 1$.  Then
$$\|\wght^{\gamma_1}w\|_{L^q(\st)}\le C_{q,\gamma}\|\wght^{\gamma_2}F
\|_{L^{q/(q-1)}(\st)},
\tag1.7$$
provided that $2\le q\le 2(n+1)/(n-1)$ and
$$\gamma_1<n(1/2-1/q)-1/2, \, \, \text{and } \, \gamma_2>1/q.
\tag1.8$$
\endproclaim

As we said earlier, one should think of (1.7) as a weighted version of estimates
of Strichartz \cite{23} for (1.6):
$$\|w\|_{L^{2(n+1)/(n-1)}(\st)}\le C\|F\|_{L^{2(n+1)/(n+3)}(\st)}.
\tag1.9$$
If one interpolates between this inequality and (1.7) one finds that the latter
holds for a larger range of weights (see also our remarks for the radial
case below).  However, for the sake of simplicity, we have only stated the
ones that will be used in our proof of Theorem 1.1.

\bigskip

Having stated our main results, let us now give the simple argument showing
how they imply Theorem 1.1.  To do so let us first notice that by shifting
the time variable by $R>0$ they yield 
$$\multline
\| \, ((t+R)^2-|x|^2)^{\gamma_1}w\, \|_{L^q(\st)}\le C
\|\, ((t+R)^2-|x|^2)^{\gamma_2}F\|_{L^{q/(q-1)}(\st)},
\\
\text{if } \, F(t,x)=0, \, \, |x|\ge t+R-1,
\endmultline
\tag1.7'$$
where $q$ and the $\gamma_j$ are as in (1.7).

It is more convenient to use this equivalent version of (1.7) in proving 
Theorem 1.1.  The key step will be to use it to establish the following

\proclaim{Lemma 1.3}  Let $u_{-1}\equiv 0$, and for $m=0,1,2,3,\dots$ let $u_m$
be defined recursively by requiring
$$\cases
\square u_m=F_p(u_{m-1})
\\
u_m(0,x)=\epsilon f(x), \, \, \partial_tu_m(0,x)=\epsilon g(x),
\endcases
$$
where $f,g\in C^\infty_0(\Rn)$ vanishing outside the ball of radius $R-1$
centered at the origin are fixed.  Then if $p_c<p\le (n+3)/(n-1)$, fix
$\gamma$ satisfying
$$1/p(p+1)<\gamma<((n-1)p-(n+1))/2(p+1)$$
and set
$$\align
A_m&=\|((t+R)^2-|x|^2)^\gamma u_m\|_{L^{p+1}(\st)}
\\
B_m&=\|((t+R)^2-|x|^2)^\gamma (u_m-u_{m-1})\|_{L^{p+1}(\st)}.
\endalign$$
Then there is an $\epsilon_0>0$, depending on $p$ $F_p$, $\gamma$ and the
data $(f,g)$ so that for $m=0,1,2,\dots$
$$A_m\le 2A_0 \, \, \text{and } \, 2B_{m+1}\le B_m, \, \, 
\text{if } \, \epsilon<\epsilon_0.
\tag1.10$$
\endproclaim

\demo{Proof}  Because of the support assumptions on the data, domain of dependence
considerations imply that $u_m$, and hence $F_p(u_m)$, must vanish if $|x|>t+R-1$.
It is also standard that the solution $u_0$ of the free wave equation $\square u_0=0$
with the above data satisfies $u_0=O(\epsilon (1+t)^{-(n-1)/2}(1+|t-|x||)^{-(n-1)/2})$.
Using this one finds that 
$$A_0\le C_0\epsilon,$$
for some uniform constant $C_0$.

To complete the induction argument let us first notice that for $j,m\ge0$, 
$u_{m+1}-u_{j+1}$ has zero Cauchy data at $t=0$ and
$\square(u_{m+1}-u_{j+1})=V_p(u_m,u_j)(u_m-u_j)$, where by (1.1),
$$V_p(u_m,u_j)=O((|u_m|+|u_j|)^{p-1}).$$
Since we are assuming that
$$\gamma<n(1/2-1/q)-1/2, \, \, \text{and } \, p\gamma>1/q, \, \,
q=p+1,$$
if we apply $(1.7')$ and H\"older's inequality we therefore obtain
$$\align
&\|((t+R)^2-|x|^2)^\gamma (u_{m+1}-u_{j+1})\|_{L^{p+1}}
\\
&\qquad
\le C_1\|((t+R)^2-|x|^2)^{p\gamma}V_p(u_m,u_j)(u_m-u_j)\|_{L^{(p+1)/p}}
\\
&\qquad
\le C_1\Bigl(C_2(\, \|((t+R)^2-|x|^2)^\gamma u_m\|_{L^{p+1}}
+\|((t+R)^2-|x|^2)^\gamma u_j\|_{L^{p+1}}\, )\Bigr)^{p-1}
\\
&\qquad\qquad\qquad\qquad
\times \, \|((t+R)^2-|x|^2)^\gamma (u_m-u_j)\|_{L^{p+1}},
\endalign$$
for certain constants $C_j$ which are uniform if above $p$, $\gamma$ and
$F_p$ are fixed.  Based on this we conclude that
$$\multline
\|((t+R)^2-|x|^2)^\gamma (u_{m+1}-u_{j+1})\|_{L^{p+1}}
\\
\le C_1\, (C_2(A_m+A_j))^{p-1} \,
\|((t+R)^2-|x|^2)^\gamma (u_{m}-u_{j})\|_{L^{p+1}}.
\endmultline
\tag1.11$$
If $j=-1$, then $A_j=0$ and hence we conclude that
$$A_{m+1}\le A_0+A_m/2 \quad \text{if } \, \,
C_1(C_2A_m)^{p-1}\le 1/2.$$
By the earlier bound for $A_0$, this yields the first part of (1.10) if
$C_1(2C_2C_0\epsilon_0)^{p-1}<1/2$.  If we take $j=m-1$ in (1.11), we also
obtain the other half of (1.10) if this condition is satisfied, which
completes the proof.  \qed
\enddemo

Using the lemma we easily get the existence part of Theorem 1.1.  If $\epsilon>0$
in (1.2) is small and if $u_m$ are as above we notice from the second half of 
(1.10) that $u_m$ converges to a limit $u$ in $L^{p+1}$ and hence in the
sense of distributions.  Since (1.1) and the bounds for $B_{m+1}$ yield
$$\|F_p(u_{m+1})-F_p(u_m)\|_{L^{(p+1)/p}}=O(2^{-m}),$$
and hence $F_p(u_m)\to F_p(u)$ in $L^{(p+1)/p}$, we conclude that $u$ must
converge to a weak solution of (1.2) which must satisfy (1.4) by the bounds
for $A_m$.  Since the proof of the bound for $B_{m+1}$ yields the 
uniqueness part, this completes our argument showing that the weighted
Strichartz estimates imply Theorem 1.1.

The rest of the paper will be concerned with the proof of Theorem 1.1.  We first
notice, after applying Stein's analytic interpolation theorem \cite{20}, that
to prove (1.7) it suffices to establish the bounds in the two extreme cases
where $q=2(n+1)/(n-1)$ or $q=2$.  Specifically, under our assumption that
$F(t,x)=0$ when $t-|x|<1$, we must show that for $n\ge2$
$$\multline
\|\wght^{\gamma_1}w\|_{L^{2(n+1)/(n-1)}(\st)}
\le C_\gamma \|\wght^{\gamma_2}F\|_{L^{2(n+1)/(n+3)}(\st)},
\\
\text{if } \, \gamma_1<(n-1)/2(n+1)<\gamma_2,
\endmultline
\tag1.12$$
and that
$$\|\wght^{-\gamma}w\|_{L^2(\st)}\le C_\gamma\|\wght^{\gamma}F\|_{L^2(\st)}, \, \,
\text{if } \, \gamma>1/2.
\tag1.13$$

Most of the rest of the paper will be devoted to the proof of (1.12).  The
$L^2$-estimate is much easier, following essentially from a twofold application
of the Sobolev trace theorem.

In proving the weighted Strichartz inequality (1.12) we shall of course exploit
our support assumption and the favorable condition on the weights.  Indeed
since $t^2-|x|^2\ge t$ on the supports of $w$ and $F$, we shall see right away
that it suffices to prove a variant of (1.12) where we assume in the left that
the norm is taken over a dyadic strip where $T/2\le t\le T$ for  $T$ large.
Assuming this, our estimate naturally splits into two pieces.  The easiest
half involves estimating the contribution to $w$ of the part of $F(t,x)$ where,
say, $t\ge T/10$.  Here, using elementary geometry and exploiting the Lorentz-invariance of
the weights, it turns out that we can reduce matters to an estimate which follows
from the usual $L^2$-calculus of Fourier integral operators.  The analysis of the
relatively small-time contributions of $F$, though, is harder since the resulting
Fourier integral operators that arise become increasingly degenerate in places as
$T\to +\infty$ and hence, as in the preceding case, we cannot hope to appeal to
H\"ormander's $L^2$-theorem.  Fortunately, though, these sorts of degenerate
Fourier integral operators have been studied before, for instance 
in Sogge and Stein \cite{19}, and
the weights in the inequalities that arise compensate for the degeneracy of the
operators, allowing the estimates to hold.  It turns out, though, that the techniques
from \cite{19} can only be used to handle the high-frequency parts of the Fourier
integrals that arise.  This in part accounts for the fact that the second step in
the proof of (1.12) is much harder than the first.  Fortunately, though, we can
handle the low frequency part using stationary phase and elementary geometric facts
which are somewhat similar to the ones mentioned before.  The two 
geometrical facts that we use, which are based on properties of the intersection
of essentially externally and internally tangent spheres,  
have widely been used in harmonic analysis, especially in 
the study of circular maximal inequalities (see \cite{1}, \cite{16}, \cite{25}).

Before turning to the details, we thought it might be well to see how under the
assumption of spherical symmetry it is easy to prove Theorem 1.2.  It turns out
that under this assumption we can also prove a stronger estimate which probably
involves the optimal range of weights.  For brevity, we shall only consider odd
spatial dimensions for the radial case.  The argument for even $n$ is a bit more
technical, due to the lack of strong Huygen's principal; however, using techniques
from \cite{9} one could adapt the proof to handle even $n$.

With this in mind, let us close this section with the following

\proclaim{Theorem 1.4} Let $n$ be odd and assume that 
$F$ is spherically symmetric and
supported in the forward light cone $\{(t,x)\in \Bbb R^{1+n}:\, |x|\le t\}$.
Then if $w$ solves $(1.6)$ and if $2<q\le 2(n+1)/(n-1)$
$$\multline
\|(t^2-|x|^2)^{-\alpha} w\|_{L^{q}(\Bbb R^{1+n}_+)} \le 
C_\gamma\|(t^2-|x|^2)^{\beta}F\|_{L^{q/(q-1)}(\Bbb R^{1+n}_+)},
\\
\text{if } \, 
\beta<1/q,\quad 
\alpha+\beta+\gamma=2/q,\quad\text{where} \quad 
\gamma=(n-1)(1/2-1/q).  
\endmultline
\tag1.14
$$
\endproclaim

\demo{Proof} 
For odd $n$ we have the formula
$$
w(t,r)
=\frac{1}{r^{(n-1)/2}} 
\int_0^t\int_{|t-r-s|}^{t+r-s}\,\, P_m(\mu)F(s,\rho)\,\rho^{(n-1)/2} d\rho ds,
$$
where $P_m(\mu)$ are Legendre polynomials of degree $m=(n-3)/2$
and $\mu=(r^2+\rho^2-(t-s)^2)/2r\rho$ satisfies $-1\leq \mu\leq 1$ in the 
domain of integration. 
(See e.g. $(3.2')$ and the formula after (3.11) in \cite{9}.) 
Multiplying by $K(t,r) (t^2-r^2)^{-\alpha}$ and integrating with respect to 
$dxdt=c_n r^{n-1} dr dt$, 
we see that we must show that 
$$
\int_{0}^{\infty}\int_0^\infty \int_0^t\int_{|t-r-s|}^{t+r-s} 
\frac{\,\,| K(t,r)| r^{(n-1)/p}  (s^2-\rho^2)^{\beta}|F(s,\rho)| \rho^{(n-1)/p}}
{(r\rho)^{\gamma}(s^2-\rho^2)^{\beta} (t^2-r^2)^{\alpha}  }\,d\rho \,ds\,  dr\, dt ,
$$
is bounded by a constant times $\| K\|_{L^{q/(q-1)}}\cdot 
\| \wght^\beta F\|_{L^{q/(q-1)}}$, if
$\gamma=(n-1)/2-(n-1)/q$ and the norms are with respect to $dxdt=c_n r^{n-1}\, dr dt$.
To do this it is convenient to introduce $u=t+r$, $v=t-r$,
$\xi=s+\rho$ and  $\eta=s-\rho$ as new variables and let
$G(\xi,\eta)=  (s^2-\rho^2)^{\beta}| F(s,\rho)|\rho^{(n-1)/p}$
and $H(u,v)= | K(t,r) | r^{(n-1)/p}$, $p=q/(q-1)$.  
We then must show that 
$$\multline
\iiiint_{0\leq \eta\leq v\leq \xi\leq u}
{\frac{G(\xi,\eta)H(u,v)\, }{|u-v|^\gamma |\xi-\eta|^\gamma |\xi\eta|^\beta |uv|^\alpha}
d\xi\, d\eta \, du\, dv }
\\
\leq C \| G\|_{L^{q/(q-1)}}\|H\|_{L^{q/(q-1)} }. 
\endmultline
\tag1.15
$$
In the domain of integration the kernel is bounded by 
$$
\frac{1}{|u-\xi|^{\gamma} |\xi|^\beta |u|^\alpha }\cdot
 \frac{1}{|v-\eta|^{\gamma} |\eta|^\beta |v|^\alpha }
$$
and (1.15) now follows from two applications of the inequality
$$
\| f\|_{L^q[0,\infty]}\leq C\|g\|_{L^{p}[0,\infty]} ,
\quad\text{if}\quad f(u)=\int_0^{u}{\frac{g(\xi)\, d\xi}{|u-\xi|^{\gamma} |\xi|^\beta |u|^\alpha}\,},
\tag1.16$$
where 
$$1<p<q<\infty, \, \,  \alpha+\beta+\gamma =1-(1/p-1/q), \, \,  \alpha+\beta\geq 0,
\, \, 
\text{and } \,  \alpha+\gamma>1/q.
$$
Notice that, for dual exponents $q$ and $p=q/(q-1)$,
$\alpha+\beta+\gamma=2/q$.
Therefore, $\alpha+\beta\geq 0$ is equivalent to 
$\gamma=(n-1)(1/2-1/q)\leq 2/q$ which holds if and only if
$q\leq 2(n+1)/(n-1)$.
In proving (1.16) we may assume that $g(\xi)\geq 0$.
Since $\alpha+\beta\geq 0$ we have  $f(u)\leq Cf_1(u)+Cf_2(u)$ where 
$$
f_1(u)=\frac{1}{|u|^{\gamma+\alpha}}\int_{0}^{u/2}{\frac{g(\xi)\, d\xi}{|\xi|^{\beta} } },
\quad f_2(u)=\int_{-\infty}^{+\infty}
{\frac{g(\xi)\, d\xi}{|u-\xi|^{\gamma+\alpha+\beta}  }}. 
$$
That  $\|f_2\|_{L^q}\leq C\|g\|_{L^{q/(q-1)}}$ is just Hardy-Littlewood's inequality for
fractional integrals. 
Dividing the integral $f_1(u)$ further into $0\leq \xi\leq u/4$ and
$u/4\leq \xi\leq u/2$, we see that
$f_1(u)\leq 2^{-(\alpha+\gamma)} f_1(u/2)+C f_2(u)$ and hence 
$$\|f_1\|_{L^q}\leq 2^{1/q-(\alpha+\gamma)} \|f_1\|_{L^q}
+C^\prime\| g\|_{L^{q/(q-1)}}.$$
Now $1/q-(\alpha+\gamma)<0$, by assumption,
so this gives the desired {\it a priori} inequality for $f_1$
and hence for $f$. 
Clearly,
$f\in L^q$, when $\alpha+\gamma>1/q$,
 if $g$ is bounded and compactly supported,  so (1.16) follows. 
\qed\enddemo

As a side remark, we note that we can use (1.14) to give an elementary proof of
John's existence theorem for $n=3$.  Indeed since the mapping from $F$ to $w$
is a positive operator when $n=3$, (1.14) yields
$$\multline
\bigl\|\, \wght^{-\alpha}\sup_{\theta\in S^2}|w(t,r\theta)|\, \bigr\|_{L^q(r^2drdt)}
\\
\le C\bigl\|\, \wght^{\beta}\sup_{\theta\in S^2}|F(t,r\theta)|\, \bigr\|_{L^{q/(q-1)}
(r^2drdt)},
\endmultline$$
for $2<q\le 4$ and $\alpha$ and $\beta$ as in (1.14).  Since this is stronger
than the estimates employed in the proof of Lemma 1.3 for $n=3$, we conclude that
in this case one always has global small-amplitude solutions of (1.2) if
$p>1+\sqrt2$.

The authors would like to thank S. Klainerman for his support and encouragement
throughout this project.

\head {\bf 2. Lorentz invariance and bounds for relatively small
times}\endhead

In proving our weighted Strichartz inequality $(1.12)$, we shall see that,
because the weights in the left are smaller than those in the right,
we can easily reduce
matters to proving estimates where in the left the norms are taken over
sets where $t$ and $t-|x|$ belong to dyadic intervals.  Let us first handle
the case where $T/2\le t\le T$, for some $T\ge2$, and $(t,x)$ belongs to
the ``middle part'' of the light cone, that is, $|x|\le t/2$.  This is the
model case.  It turns out to be the easiest to handle, and, using Lorentz
rotations as in \cite{10}, we shall reduce much of our task to this one.
Unfortunately, as we shall see, part of the weighted estimate cannot be 
handled in this manner.  However, in the next section we shall show that
the remaining cases can be handled using estimates for degenerate Fourier
integrals in the spirit of \cite{19}.

With this in mind, our first task then is to establish the following result,
which, among other things, ensures that the variant of $(1.12)$ holds where the
norm in the left is taken over all $(t,x)$ with $|x|\le t/2$.

\proclaim{Proposition 2.1}  Let $n\ge2$ and $q=2(n+1)/(n-1)$, and assume
that $F(t,x)=0$ if $t^2-|x|^2\le1$.  Then if $w$ is the solution of the 
inhomogeneous wave equation $\square w = F$ in $\st$
with zero Cauchy data
at $t=0$, 
$$\multline
\|\wght^{1/q}w\|_{L^q(\{(t,x): \, |x|\le t/2, \, T/2\le t\le T\})}
\\
\le C(\log T)^{1/q}\, \|\whgt^{1/q}F\|_{L^{q/(q-1)}}, \, \, \,
T\ge 2, 
\endmultline
\tag2.1$$
where $C$ depends only on the dimension.
\endproclaim

\demo{Proof}
Let $w_T(t,x)=w(Tt,Tx)$ and $F_T(t,x)=T^2F(Tt,Tx)$, so that 
$\square w_T=F_T$.  Then the first step is to notice that (2.1)
is equivalent to
$$
\|w_T\|_{L^q(\{(t,x): \, |x|\le t/2, \, 1/2\le t\le 1\})}
\le C(\log T)^{1/q} \|\wght^{1/q}F_T\|_{L^{q/(q-1)}}.
\tag2.1'$$
Note that $F_T=0$ if $t^2-|x|^2\le 1/T^2$.
Taking into account
the domain of dependence, we may also assume that $F_T(t,x)=0$ if
$t<1/4$ if the spatial dimension $n$ is odd.  It is not difficult
to make a similar reduction in even spatial dimensions.  
To see this, we need to recall that in any dimension 
$w_T=E_+*F_T$, where $E_+(t,x)=\pi^{(1-n)/2}/2\cdot \chi^{-(n-1)/2}_+\wght$, 
if $t\ge 0$ and $0$ otherwise.
\footnote{Here $\chi^z\wght$ denotes the pullback of the distribution
$(\Gamma(z))^{-1}x_+^z$ via the Lorentz form $t^2-|x|^2$.}
Because of this, we can
assume that $F_T$ vanishes when $t<1/8$ if we use H\"older's inequality,
since if $1/2\le t\le 1$ and $|x|\le t/2$
$$\|E_+(t-s,x-y) \, (s^2-|y|^2)^{-1/q}\|_{L^q(\{(s,y): \, s^2-|y|^2\ge 1/T^2, \, 
\, 1/T\le s\le 1/8\})} \le C(\log T)^{1/q}.$$
To prove this one just uses the fact that the $E_+$ term is bounded because
of our assumptions on $(t,x)$ and $(s,y)$.

Because of these considerations, we conclude that in proving $(2.1')$ it suffices
to assume that $F_T$ vanishes if $t^2-|x|^2\le 1/T^2$ or $t\le 1/8$.  The difficulty
then occurs because of the fact that the weights on the right side of the
inequality are small if $(t,x)$ is near the null cone.  Indeed, if, say,
$t-|x|\ge 1/8$ on the support of $F_T$, then the estimate follows from
the well known unweighted version of Strichartz \cite{23}.  Thus, we can
further assume in proving (2.1) that
$$F_T(t,x)=0\, \, \, \text{if } \, \,
t\le 1/8, \, \, \text{or } \, \, t^2-|x|^2\le 1/T^2, \, \, 
\text{or } \, \, t-|x|\ge 1/8.
\tag2.2$$
We have made this last assumption to ensure that $t-s$ is bounded from below
when $t\ge s$, $(s,y)\in \supp F_T$, and $(t,x)$ is as in the left side
of $(2.1')$.

Having set things up, we are finally ready to prove the main part of our estimate.
Recall that, if $\Hat F(s,\xi)$ denotes the spatial Fourier transform, then
$$w_T(t,x)=(2\pi)^{-n}\int_0^t\int_{\Rn}
e^{ix\cdot\xi}|\xi|^{-1}\sin((t-s)|\xi|) \Hat F_T(s,\xi)\, d\xi ds.$$
Therefore, if we let
$$(W^zF_T)(t,x)=(z-(n+1)/2)e^{z^2}\int_0^t\int_{\Rn}
e^{ix\cdot\xi}|\xi|^{-z}\sin((t-s)|\xi|)\Hat F_T(s,\xi)\, d\xi ds,$$
by complex interpolation, it suffices to show that
$$\|W^zF_T\|_{L^\infty(|x|\le t/2, \, 1/2\le t\le 1)}
\le C\|F_T\|_{L^1}, \, \, \,
\text{Re }z=(n+1)/2,
\tag2.3$$
and
$$\|W^zF_T\|_{L^2(|x|\le t/2, \, 1/2\le t\le 1)}\le C\, (\log T)^{1/2}
\|(t^2-|x|^2)^{1/2}F_T\|_{L^2}, \, \, \, \text{Re }z=0.
\tag2.4$$

Since $t-s$ is bounded from below, because of our assumptions, (2.3) follows
from the well known stationary phase estimate
$$\Bigl| \, ye^{-y^2} \int e^{ix\cdot\xi+it|\xi|}|\xi|^{-(n+1)/2+iy}\,
d\xi\, \bigr|\le C_n t^{-(n-1)/2} \, .
\tag2.3'
$$

To prove the $L^2$ estimates we note that $W^z=(W^z_++W^z_-)/2i$, where
$$\align
&(W^z_\pm F)(t,x)
\\
&\, \, \, =
(z-(n+1)/2)e^{z^2}
\iiint 
e^{i(x-y)\cdot\xi\pm i(t-s)|\xi|}F_T(s,y)|\xi|^{-z} d\xi dy ds
\\
&\, \, \, =(z-(n+1)/2)e^{z^2}\int_{1/T}^{1/8}\iint 
e^{i(x-y)\cdot\xi\pm i(t-\tau-|y|)|\xi|}F_T(|y|+\tau,y)|\xi|^{-z} d\xi dy d\tau.
\endalign$$
Here we are assuming that $(t,x)$ is as in the left side of (2.4)
so that $s$ is smaller than $t$ in the support of the first integrand.
Note that, by H\"older's inequality, the last quantity is dominated by
$(\log T)^{1/2}$ times
$$
\Bigl( 
\int \Bigl| \tau^{1/2}(z-(n+1)/2)e^{z^2}\iint
e^{i(x-y)\cdot\xi\pm i(t-\tau-|y|)|\xi|}F_T(|y|+\tau,y)|\xi|^{-z} d\xi dy
 \Bigr|^{2}d\tau \Bigr)^{1/2}.$$
Since $\tau=s-|y|$ when $s=\tau+|y|$, we conclude that it suffices to show
that for $\tau\le 1/8$ and $1/2\le t\le 1$ we have the uniform bounds
$$\multline
\|\tilde W^z_\pm f(t-\tau,\, \cdot \, )\|_{L^2(\{x: \, |x|\le t/2\})}
\\
\le C\|f\|_{L^2}, \, \, \text{Re }z=0, \, \, \supp f\subset
\{y: \, 1/8\le |y|<t-\tau\},
\endmultline
\tag2.4'$$
if
$$
(\tilde W^z_\pm f)(t-\tau,x)=(z-(n+1)/2)e^{z^2}\iint e^{i(x-y)\cdot \xi
\pm i(t-\tau-|y|)|\xi|}f(y)|\xi|^{-z} \, d\xi dy.$$

We should emphasize that this estimate would not hold if in the left the norm
were taken over all of $x\in \Rn$.  Because of our localization, though, the
bound follows from H\"ormander's theorem \cite{4} regarding $L^2$ bounds for Fourier
integrals since the symbols involved belong to a bounded subset of zero-order
symbols and since the operator has a canonical relation which is a canonical
graph in $T^*\Rn \times T^*\Rn$.  Indeed, if
$$\varphi_\pm =(x-y)\cdot\xi \pm (t-\tau-|y|)|\xi|,$$
is the phase, the last condition is equivalent to the statement that
for $\xi\ne0$
$$\text{det }\partial^2\varphi_\pm/\partial y_j\partial \xi_k\ne 0 
\, \, \text{and } \, \nabla_y\varphi_\pm\ne 0,
\, \, \, \text{if } \, \nabla_\xi \varphi_\pm =0.$$
(See, e.g. \cite{17, p. 174}.)  However, since this Hessian determinant
is just $-1\mp\langle y/|y|,\xi/|\xi|\rangle$ and since $\nabla_y\varphi_\pm
=-\xi\mp |\xi|\cdot y/|y|$,
this condition is met
since $\nabla_\xi\varphi_\pm$ does not vanish in a conic neighborhood
of $\mp y/|y|$ if $|x|\le t/2$.  For instance, if $x=0$ one must have
$\xi/|\xi|=\pm y/|y|$ if the $\xi$-gradient vanishes since, by assumption,
$t-\tau-|y|>0$.

Since we have argued that the remaining estimate $(2.4')$ follows from
the usual $L^2$ Fourier integral calculus, the proof is complete.  \qed
\enddemo

Let us now see that we can use (2.1) to estimate $w$ if the norm is
taken over a set where $T/2\le t\le T$ and $F(t,x)$ vanishes when $t$ is
smaller than a fixed multiple of $T$, if, as in Theorem 1.2, we
also assume that $|x|<t-1$ in the support of $F$.
To be more specific, if we let
$w=w^1+w^0$, where $\square w^1=F^1$ with zero data and if
$F^1(t,x)=F(t,x)$ for $t\ge T/10$, but zero otherwise then we claim
that, for $q=2(n+1)/(n-1)$,
$$\multline
\|\wght^{1/q-\epsilon}w^1\|_{L^q(\{(t,x): \, T/2\le t\le T\})}
\\
\le CT^{-2\epsilon}(\log
T)^{2/q}\|\wght^{1/q+\epsilon}F^1\|_{L^{q/(q-1)}}.
\endmultline
\tag2.5$$
Note that $w^1$ and $F^1$, like $w$ and $F$, vanish when $t-|x|\le
1$.

The next step is to also break things up with respect to the $t-|x|$
variable.  Specifically, we note that (2.5) follows from the further
localized bounds
$$\multline
\|\wght^{1/q}w^1\|_{L^q(\{(t,x): \, T/2\le t\le T, \, 2^{k-1}\le
t-|x|\le 2^k\})}
\\
\le C(\log T)^{1/q}\|\wght^{1/q}F^1\|_{L^{q/(q-1)}}.
\endmultline
\tag2.5'$$
Clearly in what follows we may assume that $2^k\le 4T$, since otherwise
the condition in the left will not be satisfied.  Also, if we set 
$T_k=T/2^k$ and let $w^1_k(t,x)=w^1(2^kt,2^kx)$ and $F^1_k(t,x)=2^{2k}
F^1(2^kt,2^kx)$, then our task is equivalent to showing that
$$\multline
\|\wght^{1/q}w^1_k\|_{L^q(\{(t,x): \, T_k/2\le t\le T_k, \, \,
1/2\le t-|x|\le 1\})}
\\ 
\le C(\log T)^{1/q}\|\wght^{1/q}F^1_k\|_{L^{q/(q-1)}}.
\endmultline
\tag2.5''$$
Note that $\wght^{1/2}\ge 2^{-k}$ and $t\ge T_k/10$ on the support of
$F^1_k$.

To use all of this we shall need the following two lemmas.

\proclaim{Lemma 2.2}  Let $E_+$ be the forward fundamental solution for
$\square$. If $0\le t-|x|\le 1$, $t/10\le s\le t$, 
and $s-1\le |y|\le s$,
then
$$|\, x/|x|-y/|y|\, |\le C/\sqrt{t}\,  \, \, \text{if } \,
(t,x,s,y)\in \supp E_+(t-s,x-y),$$
for some uniform constant $C$.
\endproclaim

\proclaim{Lemma 2.3}  Suppose that $K(x,y)$ is a measurable function on
$\Bbb R^m\times \Bbb R^n$ and set 
$$Tf(x)=\int K(x,y)f(y)\, dy.$$
Suppose further that we can write $\Bbb R^m$ and $\Bbb R^n$ as disjoint
unions $\Bbb R^m=\cup_{j\in \Bbb Z^d}A_j$ and $\Bbb R^n=\cup_{k\in \Bbb
Z^d} B_k$, where if $x\in A_j$, then $K(x,y)=0$ when $y\in B_k$ with
$|j-k|\ge C$, for some uniform constant $C$.  Then, if we let $T_{jk}$
denote the integral operator with kernel $K_{jk}$, where
$K_{jk}(x,y)=K(x,y)$
if $(x,y)\in A_j\times B_k$ and zero otherwise,
$$\|T\|_{L^p\to L^q}\le (2C+1)^d \cdot\sup_{j,k}\|T_{jk}\|_{L^p\to L^q},
$$
provided that $1\le p\le q\le \infty$.
\endproclaim

Using these two lemmas it is easy to obtain $(2.5'')$ from (2.1).
We first notice that it is enough to prove the variant of $(2.5'')$
where in the left we also assume that $|x/|x|-\nu|\le C/\sqrt{T_k}$ for
some $\nu\in S^{n-1}$.  Next, we let $\omega=(t,x)/\sqrt{t^2-|x|^2}$
denote the projection of $(t,x)$ onto the unit hyperboloid $\Bbb H^n$,
we notice that if $(t_j,x_j)$, $j=1,2$ are two points in the set where
$T_k/2\le t\le T_k$, $1/2\le t-|x|\le 1$, $|x/|x|-\nu|\le
C/\sqrt{T_k}$,
then we must have $\text{dist}(\omega_1,\omega_2)\le C_0$, for some
uniform constant with $\dist$ denoting the distance on $\Bbb H^n$ with
respect
to the restriction of the Lorentz metric $dx^2-dt^2$ to the
hyperboloid.  Hence, after making a Lorentz rotation which sends this
set
to the ``middle'' of the light cone, we see that the remaining estimate
would follow from
$$\multline
\|\wght^{1/q}w\|_{L^q(\{(t,x): \, |x|\le t/2, \, \, T_k^{1/2}/2\le t\le
T^{1/2}_k\})}
\\
\le C(\log T)^{1/q}\|\wght^{1/q}F\|_{L^{q/(q-1)}},
\endmultline$$
if $\square w=F$ with zero data and $F(t,x)=0$ if $\wght^{1/2}\le
2^{-k}$, as before.  This in turn follows from (2.1) if we rescale
since
$T_k=T/2^k$ and $2^k\le 4T$.

Thus our proof of (2.5) will be complete once we have established the
above elementary lemmas.

The first one is quite standard and relies on a geometric fact that has
been used extensively in the study of circular maximal operators and 
related topics.  See, e.g., \cite{1}, \cite{16} and \cite{25}.

\demo{Proof of Lemma 2.2}  
The conclusion trivially holds for a large constant $C$ if $t$ is small,
so in what follows we shall assume, say, $t\ge20$, so that our assumptions
then give $2|y|\ge s$.
We then need to use the following version of
Huygen's principle:
$$E_+(t-s,x-y)=0, \quad \text{if } \, |x-y|> t-s.$$
Using the identity
$$
|x-y|^2=(|x|-|y|)^2+2(|x|\, |y|-x\cdot y)=
(|x|-|y|)^2+|x| \, |y|\, |\, x/|x|-y/|y|\, |^2,
$$
we see that $|x-y|^2\le (t-s)^2$
is equivalent to
$$\Bigl|\, \frac{x}{|x|}-\frac{y}{|y|}\Bigr|^2
\le \frac{(t-s)^2-(|x|-|y|)^2}{|x|\, |y|}
=\frac{(t-|x|-(s-|y|))(t+|x|-(s+|y|))}{|x|\, |y|}.$$
Since $|y|\le s$ the right side is $\le (t-|x|)(t+|x|)/|x|\, |y|$,
which in turn is $O(t^{-1})$ if the assumptions are fulfilled.  \qed
\enddemo

Notice how the lower bound for $s$ is essential.  It is for this reason
that we must use different techniques to estimate the norm of $w$ over
$T/2\le t\le T$ if $F$ is supported in a region where $t$ is much
smaller than $T$.

We still must handle the last lemma:

\demo{Proof of Lemma 2.3}  Let us assume that $q<\infty$, since the proof for
$q=\infty$
is similar.  We first notice that H\"older's inequality gives
$$\align
\int |Tf(x)|^q\, dx \, &=\, \sum_j \int_{A_j}|Tf(x)|^q\, dx
\\
&\le (2C+1)^{(q-1)d} \, \sum_{\{(j,k): \, |j-k|\le C\}}
\int |T_{jk}f_k(x)|^q\, dx,
\endalign$$
where $f_k(y)=f(y)$ if $y\in B_k$ and zero otherwise.  This in turn is
$$\align
&\le (2C+1)^{(q-1)d}\,\sup \|T_{jk}\|^q_{L^p\to L^q} \cdot 
\sum_{\{(j,k): \, |j-k|\le C\}}\bigl(\, \int |f_k(y)|^p\, dy\,
\bigr)^{q/p}
\\
&\le (2C+1)^{qd}\, \sup \|T_{jk}\|_{L^p\to L^q}^q \cdot \sum_{k\in \Bbb
Z^d} \|f_k\|_{L^p}^q
\\
&\le (2C+1)^{qd}\, \sup \|T_{jk}\|_{L^p\to L^q}^q \cdot 
\|f\|_{L^p}^q\, ,
\endalign$$
using our assumption that $p\le q$ in the last step.
\qed
\enddemo

\head{\bf 3. Degenerate Fourier integrals and bounds for relatively large
times}\endhead

To finish our proof of $(1.12)$ we have to estimate $w_0$ which involves the
contributions to $w$ from relatively small-time parts of $F$.  Specifically,
if $T\ge10$, and if we set $F^0(t,x)=F(t,x)$ if $t\le T/10$ and $0$ otherwise
and if $w^0$ is the solution of $\square w^0=F^0$ with zero data then it suffices
to show that
$$\multline
\|\wght^{1/q-\epsilon}w^0\|_{L^q(\{(t,x): \, T/2\le t\le T\})}
\\
\le CT^{-\epsilon/4}\|\wght^{1/q+\epsilon}F^0\|_{L^{q/(q-1)}}.
\endmultline
\tag3.1$$
As before, $q=2(n+1)/(n-1)$.  Note that $F^0$ and $w^0$, like $F$ and $w$ in
$(1.12)$, vanish if $t-|x|\le 1$.  Clearly since $w=w^0+w^1$, this inequality along
with (2.5) yields $(1.12)$.

The proof of (3.1) is in many ways opposite to that of (2.5).  Instead of relying
on $L^2$ estimates for ``non-degenerate'' Fourier integrals, the main part here 
rests on weighted $L^2$ estimates for the degenerate Fourier integral operators
which arise in the study of the characteristic Cauchy problem.  Also, the main
reduction now will rely on the geometry of internally tangent spheres, rather
than externally tangent ones as in the earlier estimate.

To set up the main estimate, let us make a couple of reductions which
exploit
the fact that the weights in (3.1) scale favorably because of the
$\epsilon$ parts.  First, if we assume additionally that $F^0$
vanishes for $t\notin [T_0,2T_0]$, then it suffices to show that the
variant of (3.1) holds where $T^{-\varepsilon/4}$ is replaced by
$(TT_0)^{-\epsilon/4}$ in the right.  If we assume further that $F^0$
also vanishes if $t-|x|\notin [\delta_0T_0,2\delta_0T_0]$ then it
suffices to show that the inequality holds with operator norm
$O(T^{-\epsilon/4} T_0^{-\epsilon/2})$.  Since by domain of dependence
considerations, $w^0$ will then vanish if $t-|x|\le \delta_0T_0$, we
conclude that this in turn would follow from showing that for
$\delta\ge \delta_0$
$$\multline
\|\wght^{1/q-\epsilon}w^0\|_{L^q(\{(t,x): \, T/2\le t\le T, \, \,
\delta T_0\le t-|x|\le 2\delta T_0\})}
\\
\le C(TT_0)^{-\epsilon/2}\|\wght^{1/q+\epsilon}F^0\|_{L^{q/(q-1)}},
\endmultline
\tag3.1'$$
assuming as we are now that
$$F^0(t,x)=0\, \, \, \text{if } \, \, t\notin [T_0,2T_0], \, \, 
\text{or } \, \, t-|x|\notin [\delta_0T_0,2\delta_0T_0].$$
Note that we must have $\delta_0\ge 1/T_0$.

One advantage of this inequality is that in both sides the weights are
essentially constant on the supports.  Specifically, our task amounts
to showing that
$$\multline
(TT_0\delta)^{1/q-\epsilon}\|w^0\|_{L^q(\{(t,x): \, T/2\le t\le T, \, \,
\delta T_0\le t-|x|\le 2\delta T_0\})}
\\
\le C(TT_0)^{-\epsilon/2}\, (T_0^2\delta_0)^{1/q+\epsilon} \, 
\|F^0\|_{L^{q/(q-1)}}.
\endmultline
$$
Since $1/T_0\le \delta_0\le \delta$, by rearranging terms, this in turn
would follow from
$$\multline
(T/T_0)^{1/q-\epsilon/2} \, \delta^{1/q+\epsilon/2}
\|w^0\|_{L^q(\{(t,x): \, T/2\le t\le T, \, \,
\delta T_0\le t-|x|\le 2\delta T_0\})}
\\
\le C\delta_0^{1/q}\|F^0\|_{L^{q/(q-1)}}.
\endmultline$$
Finally, if we let $G(t,x)=T_0^2F^0(T_0t,T_0x)$ and
$v(t,x)=w^0(T_0t,T_0x)$ so that
$$\square v=G, \, \, v(0,\cd)=\partial_t v(0,\cd)=0,$$
and 
$$\supp G\subset \{(t,x): \, 1\le t\le 2, \, \, 
\delta_0\le t-|x|\le 2\delta_0\},$$
then, if we abuse notation and let $T$ now denote $T/T_0$, the last
inequality is in turn equivalent to
$$\align
&T^{1/q-\epsilon/2}\delta^{1/q+\epsilon/2}\, \|v\|_{L^q(\{(t,x): \, 
T/2\le t\le T, \, \, \delta \le t-|x|\le 2\delta \} )}
\tag3.2
\\
&\qquad\qquad\quad\le C\delta_0^{1/q}\|G\|_{L^{q/(q-1)}}.
\endalign
$$
Here we can assume that $\delta_0\le \delta$, and, since we have
replaced
$T/T_0$ by $T$, our assumption on $T$ is now that $T\ge 10$.

It is easy to handle the extreme cases of this inequality where, say,
$\delta_0\le \delta\le 10\delta_0$, or $\delta\ge 10$.

For the first case, a stronger version would say that, for $T\ge10$, 
$$T^{1/q}\|v\|_{L^q(\{(t,x): \, T/2\le t\le T\})}
\le C\|G\|_{L^{q/(q-1)}}, \, \, \text{if } \, 
G(t,x)=0, \, \, t\notin [1,2].$$
But if we use a routine freezing argument (see, e.g. \cite{17,
\S 0.3}), we see that this follows from the following estimates of 
Strichartz \cite{23}
$$\|u(t-s,\cd)\|_{L^q(\Bbb R^n)}\le
C|t-s|^{-2/q}\|g\|_{L^{q/(q-1)}(\Bbb R^n)},$$
where 
$$u(t,x)=(2\pi)^{-n}\int e^{ix\cdot\xi}\sin(t|\xi|)\Hat g(\xi) \, 
d\xi/|\xi|.$$
This inequality implies the preceding one since if we let 
$K(t,s)=|t-s|^{-2/q}$, when $(t,s)\in [T/2,T]\times [1,2]$ and $0$
otherwise then, by H\"older's inequality, the associated integral
operator sends $L^{q/(q-1)}(\Bbb R)$ to $L^q(\Bbb R)$ with norm
$O(T^{-1/q})$.

The case where, in (3.2), $\delta\ge 10$ is even easier to handle.
Indeed, since the forward fundamental solution $E_+(t,x)$
vanishes for $t<0$ and
for $t\ge0$ is a multiple of
$\chi_+^{-(n-1)/2}(t^2-|x|^2)$, a calculation shows that for $1\le s\le
2$, $|y|\le s$, and $n\ge2$,
$$\int_{T/2}^T\int_{t-|x|\ge10}\bigl|\wght^{1/q}E_+(t-s,x-y)|^q\, 
dtdx =O(1),
$$
if as above $q=2(n+1)/(n-1)$.  This just follows from the fact that
the $E_+$ term is $O((t(t-|x|))^{-(n-1)/2})$ because of our
assumptions.
If we use H\"older's inequality as in the proof of Proposition 2.1, we
conclude that, as claimed, (3.2) must hold when $t-|x|\ge 10$.

To handle the remaining cases where $10\delta_0\le \delta\le 10$, first
notice that if we use H\"older's inequality, as in the
proof of Proposition 2.1, then we find that $v$ in $(3.2)$ is dominated
by $\delta_0^{1/q}$ times
$$\Bigl(\int_{\delta_0}^{2\delta_0}\Bigl|\iint
e^{i(x-y)\cdot\xi}|\xi|^{-1}
\sin((t-\tau-|y|)|\xi|)G(\tau+|y|,y)d\xi
dy\Bigr|^{q/(q-1)}d\tau\Bigr)^{(q-1)/q}.$$
Therefore, since we are assuming that $\delta\ge 10\delta_0$, if we
replace
$t$ by $t-\tau$, we conclude that the remaining cases of (3.2) would be
a consequence of the following

\proclaim{Proposition 3.1}  For $n\ge2$ set
$$(Tg)(t,x)=\int_{\Rn}\int_{\{y\in\Rn: \, 1\le |y|\le2\}}
e^{i(x-y)\cdot\xi-i(t-|y|)|\xi|}g(y) dy d\xi/|\xi|.$$
Then, if $q=2(n+1)/(n-1)$, $\epsilon>0$,  $t>5$ and $\delta<10$
$$\|Tg(t,\cd)\|_{L^q(\{x:\, \delta\le t-|x|\le 2\delta\})}
\le C t^{\epsilon-2/q}\delta^{-\epsilon-1/q}\|g\|_{L^{q/(q-1)}}.
\tag3.3$$
\endproclaim

As before, we shall prove this using complex interpolation.  To this
end, let us set
$$
(T_zg)(t,x)
=(z-(n+1)/2)e^{z^2}\iint_{1\le |y|\le 2}e^{i(x-y)\cdot\xi-i(t-|y|)|\xi|}
|\xi|^{-z}g(y)dyd\xi,
$$
so that $T_1$ is a multiple of $T$.  Therefore, if we apply complex
interpolation we conclude that (3.3) would be a consequence of 
$$\|T_zg(t,\cd)\|_{L^\infty(\Rn)}\le Ct^{-(n-1)/2}\|g\|_{L^1}, \, \, 
\text{Re }z=(n+1)/2,
\tag3.4$$
and
$$\|T_zg(t,\cd)\|_{L^2(\{x: \, \delta\le t-|x|\le 2\delta\})}
\le C t^{\epsilon/2} \delta^{-\epsilon-1/2}\|g\|_{L^2}, \, \, 
\text{Re }z=0.
\tag3.5$$

Inequality (3.4) is a simple consequence of $(2.3')$ and our assumption
that
$t\ge5$.  The $L^2$ estimate is more delicate.  For it, we shall need
to use a bit of microlocal analysis.  These techniques will only work
for large frequencies $\xi$, depending on the scales $\delta$ and $t$.
Fortunately, it is easy to deal with the part of our operator coming
from small $\xi$ using the Sobolev trace theorem.

Let us be more specific.  To simplify the notation to follow, let us
set
$$\alpha=1+\epsilon/2.$$
If we then fix $\rho\in C^\infty$ satisfying $\rho(\xi)=0$ for
$|\xi|\le 1$ and $\rho=1$ for $|\xi|\ge2$ we claim that
$$\multline
(R_zg)(t,x)
\\
=(z-(n+1)/2)e^{z^2}\iint_{1\le |y|\le2}e^{i(x-y)\cdot\xi
-i(t-|y|)|\xi|}|\xi|^{-z}(1-\rho(t^{1-\alpha}\delta^\alpha\xi))
g(y)dy d\xi,
\endmultline$$
satisfies 
$$\|R_zg(t,\cd)\|_{L^2(\Rn)}\le Ct^{(\alpha-1)/2}\delta^{-\alpha/2}
\|g\|_{L^2(\Rn)}, \, \, \text{Re }z=0.
\tag3.6$$
But this follows by duality from the special case corresponding to
$T=2$ of the following lemma which, for future use, we state in greater
generality than is needed here.

\proclaim{Lemma 3.2}  If $T\ge 1$ then
$$\multline
\bigl\|\, \int e^{ix\cdot \xi-i|x|\, |\xi|} \Hat f(\xi) d\xi\,
\|_{L^2(T/2\le |x|\le T)} 
\\
\le CT^{1/2}\bigl( \, \|\Hat
f\|_{L^2(|\xi|\le 1)} +\sum_{k=0}^\infty 2^{k/2}\|\Hat f\|_{L^2(2^k\le
|\xi|\le 2^{k+1})} \, \bigr).
\endmultline$$
\endproclaim

\demo{Proof}  If we change variables, we can write the left side as
$$T^{-n/2}\big\|\, \int e^{ix\cdot\xi-i|x|\, |\xi|}\Hat f(\xi/T)d\xi\, 
\bigr\|_{L^2(1/2\le |x|\le 1)}.$$
If, for fixed $1/2\le t\le 1$, we apply the
the Sobolev trace theorem (see, e.g.,  \cite{5, Appendix B})
to the function $x\to \int e^{ix\cdot\xi-it|\xi|}\Hat f(\xi/T)d\xi$, we 
find that
$$\multline
\int_{\theta\in S^{n-1}}\Bigl|\int_{\Rn}e^{it\theta\cdot\xi-it|\xi|}\Hat f(\xi/T)
d\xi\Bigr|^2d\theta 
\\
\le C\Bigl(\, \|\Hat f(\xi/T)\|_{L^2(|\xi|\le1)}
+\sum_{j=0}^\infty 2^{j/2}\|\Hat f(\xi/T)\|_{L^2(2^j\le |\xi|\le 2^{j+1})}\, 
\Bigr)^2.
\endmultline$$
If we now integrate over $1/2\le t\le 1$, 
we conclude that the left side
of the inequality in the statement of the lemma is dominated by
$$\align
&T^{-n/2}\bigl(\, \|\Hat f(\xi/T)\|_{L^2(|\xi|\le 1)}+\sum_{j=0}^\infty
2^{j/2} \|\Hat f(\xi/T)\|_{L^2(2^j\le |\xi|\le 2^{j+1})}\, \bigr)
\\
&\quad\le T^{-n/2}\bigl(\, T^{1/2} \|\Hat f(\xi/T)\|_{L^2(|\xi|\le T)} 
+\sum_{2^{j-1}\ge T}2^{j/2}\|\Hat f(\xi/T)\|_{L^2(2^j\le |\xi|\le
2^{j+1})}\, \bigr)
\\
&\quad\le 2T^{1/2}\bigr( \, \|\Hat f\|_{L^2(|\xi|\le 1)}+\sum_{k=0}^\infty2^{k/2}
\|\Hat f(\xi)\|_{L^2(2^k\le |\xi|\le 2^{k+1})}\, \bigr),
\endalign$$
as desired.  \qed
\enddemo

In view of (3.6), we conclude that (3.5) would follow if
$$\|S_zg(t,\cd)\|_{L^2(\{x: \, \delta\le t-|x|\le 2\delta\})}
\le C\delta^{-1/2}\|g\|_{L^2}, \, \, \text{Re }z=0,
\tag3.7$$
if
$$(S_zg)(t,x)=e^{z^2/2}\iint_{1\le
|y|\le2}e^{i(x-y)\cdot\xi-i(t-|y|)|\xi|}
|\xi|^{-z}\rho(t^{1-\alpha}\delta^\alpha\xi)g(y) dyd\xi.$$
Note that the bounds in (3.7) are stronger than those in (3.5) or
(3.6); however, unlike in the preceding inequality, it is necessary to
assume that $t-|x|$ is larger than $\delta$ in the norm on the left.

To proceed we shall require a couple of elementary lemmas.  The first
one is the following 

\proclaim{Lemma 3.3}  If $a(\xi)$ belongs to a bounded subset of $S^0$,
and if $\rho\in C^\infty$ satisfies $\rho(\xi)=0$ for $|\xi|\le1$,  and
$\rho=1$, $|\xi|\ge2$, then, for $\alpha>1$ and $t>1$,
$$\bigl| \int e^{ix\cdot\xi
-it|\xi|}a(\xi)\rho(t^{1-\alpha}\delta^\alpha\xi)
\, d\xi \bigr| \le C_{N,\alpha} (\delta/t)^N, \, \, \,
\text{if } |\, |x|-t\, |\ge \delta/2.$$
\endproclaim

\demo{Proof}  After changing scales, we may take $t=1$.  If we then
replace $x$, $t$ and $\delta$ by $x/t$, $1$ and $\delta/t$,
respectively, it suffices to show that if $\tilde a(\xi)=a(\xi/t)$,
then
$$\bigl|\,\int e^{ix\cdot \xi-i|\xi|}\tilde a(\xi) 
\, \rho(\delta^\alpha\xi) \, d\xi\, \bigr| \le C_{N,\alpha} \delta^N,
\, \, \, \text{if } \, |\, |x|-1\, |\ge\delta/2.$$

It is easy to see, simply by integrating by parts, that these bounds
hold if, say $|x|\notin [1/2,3/2]$.  Assuming that $|x|\in[1/2,3/2]$, we
can use polar coordinates, $\xi=\lambda\theta$, $\theta\in S^{n-1}$,
and stationary phase (see e.g. \cite{17, Theorem 1.2.1}) to rewrite
our oscillatory integral as
$$\sum_\pm \int_0^\infty e^{i\lambda(1\pm|x|)}b_\pm(x,\lambda)
\lambda^{(n-1)/2}d\lambda,$$
where, because of our assumptions on the original symbol, $b_\pm=0$
for $\lambda\le C\delta^{-\alpha}$, and $(\partial/\partial\lambda)^jb_\pm
=O(\lambda^{-j})$.  Therefore, if we integrate by parts $N$ times, we see
that the preceding term is dominated for a given large $N$ by
$$\int_{C\delta^{-\alpha}}^\infty \bigl|\, 1\pm|x|\, \bigr|^{-N} \, 
\lambda^{(n-1)/2-N}\, d\lambda =O(\delta^{-N}\delta^{\alpha N
-(n+1)/2}),$$
which gives us the desired bounds since $\alpha>1$. \qed
\enddemo


To use this, let $K_z$ denote the kernel of $S_z$, that is,
$$K_z(t;x,y)=e^{z^2/2}\int e^{i(x-y)\cdot\xi-i(t-|y|)|\xi|}
|\xi|^{-z}\rho(t^{1-\alpha}\delta^\alpha\xi) d\xi.$$
We then conclude that
$$K_z=O((\delta/t)^N) \, \, \forall N, \, \, 
\text{if } \, | \, |x-y|-|t-|y|| \, |\ge\delta/2.
\tag3.8
$$
To apply this we require the following

\proclaim{Lemma 3.4}  Suppose that $t>5$, $1\le |y|\le 2$ and that 
$| \, |x-y|-|t-|y|| \, |\le \delta/2$ and $\delta\le t-|x|\le 2\delta$.
It then follows that if $\delta$ is smaller than a fixed positive
constant
$$|\, y/|y|-x/|x|\, | \in [C^{-1}_0\delta^{1/2},C_0\delta^{1/2}],$$
for some absolute constant $C_0$.
\endproclaim

The condition $|\, |x-y|-|t-|y||\, |\le \delta/2$ says
that $x$ is a distance $\le \delta/2$ from
the sphere of radius $t-|y|$ which is internally
tangent at the point $ty/|y|$ to the sphere of radius $t$ centered at
the origin.  Thus, the conclusion of the lemma is that these two
spheres
separate of distance $\approx\delta$ at points of angle
$\approx\delta^{1/2}$
from $ty/|y|$.  This type of result can also be found in \cite{1},
\cite{16} and \cite{25}.  However, for the sake of completeness,
let us give the simple proof.

\demo{Proof of Lemma 3.4}  As in the proof of Lemma 2.2 we shall use
the identity
$$\Bigl| \, \frac{x}{|x|}-\frac{y}{|y|}\, \Bigr|^2
=\frac{|x-y|^2-(|x|-|y|)^2}{|x|\, |y|}=\frac{|x-y|+|x|-|y|}{|x|\, |y|}
\cdot \frac{|x-y|-(|x|-|y|)}{\delta}\cdot \delta.$$
By our assumptions the first factor on the right is bounded from above
and below.  Writing
$|x-y|-(|x|-|y|)=|x-y|-(t-|y|)+t-|x|$, we reach the same conclusion for
the second factor, yielding the result.
\qed
\enddemo

In view of (3.8) and the overlap lemma, Lemma 2.3, we conclude 
from Lemma 3.4
that, for
small $\delta$, to prove (3.7)
it suffices to show that if $\nu\in S^{n-1}$
$$\|(S_zg)(t,\cd)\|_{L^2(\{x: \, |x/|x|-\nu|\ge\delta^{1/2}, \, \,
|x|\ge4
\})}
\le C\delta^{-1/2}\|g\|_{L^2},
\tag3.7'$$
assuming that
$$g(y)=0, \, \, 
\text{if } \, |y/|y|-\nu|\ge c_0\delta^{1/2}, \, \, \text{or } \, 
|y|\notin [1,2],$$
with $c_0>0$ being a fixed small constant.  

Our final reduction then involves the following

\proclaim{Lemma 3.5}  Suppose that $\psi(\tau)\in C^\infty(\Bbb R)$
vanishes near $\tau=0$ and equals $1$ when $|\tau|$ is large.  
Then if $\delta$ is small, $\rho$ and
$\alpha>1$ are as above, and $t>1$,
$$\Bigl|\, \int e^{ix\cdot\xi-it|\xi|}a(\xi) \psi\bigl(\delta^{-1/2}
(x/|x|-\xi/|\xi|)\bigr)\rho(t^{1-\alpha}\delta^\alpha\xi)\, d\xi \,
\Bigr|
\le C_N \, (\delta/t)^N \, ,$$
where, for a given $N$, the constants depend only on $\dist (0, \supp
\psi)$ and the size of finitely many derivatives of $\psi$, 
if $\rho$ is fixed and $a$ belongs to a bounded subset of $S^0$.
\endproclaim

\demo{Proof}  
If we let $y=x/t$ and $\tilde a(\xi)=a(\xi/t)$, the quantity we wish to
estimate can be rewritten as
$$t^{-n}\int e^{iy\cdot\xi-i|\xi|}\tilde a(\xi)\psi(\delta^{-1/2}(y/|y|-
\xi/|\xi|)) \rho((\delta/t)^\alpha\xi)d\xi.$$
We then note that
$$e^{iy\cdot\xi-i|\xi|}=-|\, y-\xi/|\xi|\, |^{-2}\Delta_\xi e^{iy\cdot\xi-i|\xi|}
-i(n-1)|\xi|^{-1}e^{iy\cdot\xi-i|\xi|}.$$
Therefore if we let
$$L(\xi,D_\xi)h(\xi)=-\Delta_\xi\bigl(\, |\, y-\xi/|\xi|\, |^{-2}h(\xi)\, \bigr)
-i(n-1)|\xi|^{-1}h(\xi),$$
the oscillatory integral we wish to estimate can be rewritten as
$$t^{-n}\int e^{iy\cdot\xi-i|\xi|}
L^N(\, \tilde a(\xi)\psi(\delta^{-1/2}(y/|y|-\xi/|\xi|))\rho((\delta/t)^\alpha\xi)\, )d\xi.$$
Note on the support of the integral $|\, y-\xi/|\xi|\, |$ is bounded below
by a uniform multiple of $\delta^{1/2}$ and so
$$\bigl|\, (\partial/\partial\xi)^\gamma(\, |\, y-\xi/|\xi|\, |^{-2})\, \bigr|
\le C_\gamma \delta^{-1-|\gamma|/2}|\xi|^{-|\gamma|}.$$
We also clearly have
$$|(\partial/\partial\xi)^\gamma(\, \tilde a(\xi)\psi(\delta^{-1/2}(x/|x|
-\xi/|\xi|))\rho((\delta/t)^\alpha\xi)\, )|\le C_\gamma\delta^{-|\gamma|/2}
|\xi|^{-|\gamma|}.$$
From this we conclude that 
$$|L^N(\, \tilde a(\xi)\psi(\delta^{-1/2}(x/|x|-\xi/|\xi|))\rho((\delta/t)^\alpha\xi)\, )|
\le C_N\delta^{-N}|\xi|^{-N},$$
which implies that for a given large $N$ the oscillatory integral is dominated
by 
$$t^{-n}\delta^{-N}\int_{|\xi|\ge C(\delta/t)^{-\alpha}}|\xi|^{-N}d\xi
=O(t^{-n-N} \cdot(\delta/t)^{(\alpha-1)N-n\alpha}),$$
yielding the desired bounds since $\alpha$ and $t$ are larger than $1$.  \qed
\enddemo


To use this lemma, note that if $1\le |y|\le2$, $|x|\ge4$,
$\nu\in S^{n-1}$ and
$|x/|x|-\nu|\ge \delta^{1/2}$ then $|(x-y)/|x-y|-\nu|\ge
\delta^{1/2}/2$
if $\delta$ is small and $|y/|y|-\nu|\le c_0\delta^{1/2}$, with
$c_0>0$ being a small uniform constant.  With this in mind, 
we conclude that, 
for small $\delta$, $(3.7')$ (and hence
(3.7)) is a consequence of the following

\proclaim{Proposition 3.6}  Suppose that 
$$f(y)=0 \, \, \, \text{if } \, \, |y|\notin [1,2] \, \, \,
\text{or } \, \, |y/|y|-e_1|\ge c_0\delta^{1/2},$$
where $e_1=(1,0,0,\dots,0)$.
Then if $c_0>0$ is smaller than a uniform
constant which is independent of $\delta<1$
$$
\int_{|\xi/|\xi|-e_1|\ge \delta^{1/2}}
\Bigl|\, \int e^{iy\cdot\xi-i|y|\, |\xi|}f(y)\, dy\, \Bigr|^2\, d\xi
\le C\delta^{-1}\, \|f\|^2_{L^2}.
\tag 3.9$$
\endproclaim

\demo{Proof}  By decomposing the conic region $\{\xi: \,
|\xi/|\xi|-e_1\, |\ge
\delta^{1/2}\}$ into a finite number of pieces, we see that it suffices
to prove the estimate when we integrate over a convex conic subset
$\Gamma_\delta$.  Note then, for later use, that there is a uniform
constant $C_1$ so that if $\delta<1$
$$|\, \zeta'/|\zeta|\, |\le C_1\delta^{-1/2}\, |\, 1-\zeta_1/|\zeta|\,
|,
\, \, \zeta\in \Gamma_\delta.
\tag3.10$$

To be able to apply an integration by parts argument we need to make
one further reduction.  Specifically, suppose that $0\le a_\delta\in
C^\infty$ is supported in the set where $1/2\le |y|\le 4$ and $|y/|y|-
e_1|\le 2c_0\delta^{1/2}$ and satisfies the natural bounds 
$$|\, (\partial/\partial y_1)^j(\partial/\partial y')^\alpha
a_\delta(y)\, | \le C_{j,\alpha} \delta^{-|\alpha|/2} \, ,
\, \, \, \forall j, \alpha,
$$
associated with this support assumption.  Here $y'=(y_2,\dots,y_n)$.
If we then set 
$$(S_\delta f)(\xi)=\int e^{iy\cdot\xi-i|y|\, |\xi|} a_\delta(y) f(y)
\, dy,$$
then it suffices to show that
$$\delta \int_{\xi\in \Gamma_\delta}|S_\delta f(\xi)|^2 \, d\xi \le C
\|f\|^2_{L^2}.
\tag 3.9'$$
The dual version of this is equivalent to 
$$\delta \|S_\delta S^*_\delta h\|_{L^2(\Gamma_\delta)}\le
C\|h\|_{L^2}, \, \, \, \supp h\subset \Gamma_\delta,
\tag 3.9''$$
where $S_\delta S^*_\delta$ is the integral operator with kernel
$$K_\delta(\xi,\eta)= \int e^{i\Phi(y,\xi,\eta)} \, a^2_\delta(y)\,
dy\, , \quad \xi,\eta\in \Gamma_\delta,$$
with the phase being
$$\Phi(y,\xi,\eta)=y\cdot (\xi-\eta)-|y|\, (\, |\xi|-|\eta|\, ).$$

Recall that $a_\delta(y) =0$ if $|y/|y|-e_1|\ge 2c_0\delta^{1/2}$.
Assuming, as we may, that $c_0$ is small enough, we claim that there is
a constant $A$ so that, for every $N$,
$$K_\delta(\xi,\eta)\le C_N
\cases 
\delta^{(n-1)/2}(1+\delta|\xi_1-\eta_1|)^{-N}, \, \, 
\text{if } \, \, \delta^{1/2}|\xi_1-\eta_1|\ge A|\xi'-\eta'|
\\
\delta^{(n-1)/2-N}(\delta^{-1}+|\xi'-\eta'|^2)^{-N}, \, \,
\text{if } \, \, \delta^{1/2}|\xi_1-\eta_1|\le A|\xi'-\eta'| \, .
\endcases
\tag 3.11$$
This yields $(3.9'')$ by Young's inequality since for large $N$
$$\multline
\delta^{(n-1)/2}\int_{A|\xi'-\eta'|\le \delta^{1/2}|\xi_1-\eta_1|}
(1+\delta|\xi_1-\eta_1|)^{-N}\, d\xi
\\ + \, 
\delta^{(n-1)/2-N}\int_{\delta^{1/2}|\xi_1-\eta_1|\le A|\xi'-\eta'|}
(\delta^{-1}+|\xi'-\eta'|^2)^{-N}\, d\xi
=O(\delta^{-1})\, .
\endmultline$$

To prove the first bound we need to integrate by parts with respect to
$y$.  To do so we note that, by the mean value theorem,
$$\align
|\partial\Phi/\partial y_1|&=\bigl| \, (\xi_1-\eta_1)-y_1/|y|\,
(|\xi|-|\eta|)\,
\bigr|
\\
&=\bigl|\, (\xi_1-\eta_1)-y_1/|y|\cdot
|\zeta|^{-1}\zeta\cdot(\xi-\eta)\,
\bigr|
\\
&\ge |\xi_1-\eta_1|\cdot |\, 1-\zeta_1/|\zeta|\, |- |\,
|\, \zeta|^{-1}\zeta'
\cdot(\xi'-\eta')\, |,
\endalign$$
where $\zeta$ is a point on the line
segment connecting $\xi$ and $\eta$.  Since we are assuming that $\Gamma_\delta$
is convex we must have $\zeta\in \Gamma_\delta$ and so
$|\, 1-\zeta_1/|\zeta|\, |\ge c\delta$ for some uniform $c>0$.
Therefore, if we let $A=2C_1$, where $C_1$ is as in (3.10), we conclude
that for $\xi, \eta\in \Gamma_\delta$ we must have
$$|\partial\Phi/\partial y_1|\ge c\delta/2\cdot |\xi_1-\eta_1|
\, \, \text{if } \, \, \delta^{1/2}|\xi_1-\eta_1|\ge A|\xi'-\eta'|.$$
Notice also that, for such $\xi$ and $\eta$,
$$|(\partial/\partial y_1)^j\Phi|\le C_j \delta |\xi-\eta|\le C_j'\delta
|\xi_1-\eta_1|, 
\, \, j\ge2, \, \, y\in \supp a_\delta.$$
If we note that
$$e^{i\Phi}=\bigl(\, 1+\bigl|\frac{\partial\Phi}{\partial y_1}\bigr|^2\, \bigr)^{-1}
\left( \bigl(1-\frac{\partial^2}{\partial y_1^2}\bigr)e^{i\Phi}
+i\frac{\partial^2\Phi}{\partial y_1^2}e^{i\Phi} \right),$$
then we can integrate by parts to see
that, for a given $N$, $K_\delta$ can be written as a combination of
terms of the form
$$\int e^{i\Phi} \, \frac{(\partial/\partial y_1)^{l_1}\Phi\dots (\partial/\partial
y_1)^{l_m}\Phi}{(1+|\partial\Phi/\partial y_1|^2)^{j+k}} \,
\bigl( \frac{\partial}{\partial y_1}\bigr)^{l_{m+1}}a^2_\delta(y)\,
dy,$$
where
$$j+k=2N, \, \, 0\le m\le j\, , \, \, \, l_j\ge 2, \, \, j\le m.$$
From this we obtain the first bounds for $K_\delta$ in $(3.11)$ since
$a_\delta$ is supported in a set of measure $O(\delta^{(n-1)/2})$.

The argument for the other bound in $(3.11)$ is similar except here we
must use our assumption that $a_\delta=0$ when $|y/|y|-e_1|\ge
2c_0\delta^{1/2}$ with $c_0$ small.  To use this, we first note that
$$\align
|\nabla_{y'}\Phi|&\ge |\xi'-\eta'|-|\, y'/|y|\, |\cdot |\, |\xi|-|\eta|\, |
\\
&\ge |\xi'-\eta'|-2c_0\delta^{1/2}|\xi-\eta|.
\endalign$$
Hence if $|\xi_1-\eta_1|\le A\delta^{-1/2}|\xi'-\eta'|$, where $A$ is
the fixed constant chosen in the last step, we conclude that
$$|\nabla_{y'}\Phi|\ge |\xi'-\eta'|/2, \, \, y\in \supp a_\delta
\tag 3.12$$
if $c_0$ is small.  Notice also that, because of our assumptions,
$$|(\partial/\partial y')^\alpha \Phi|\le C|\xi-\eta|
\le C'\delta^{-1/2}|\xi'-\eta'|.
\tag 3.13$$

To apply this, we first observe that
$$e^{i\Phi}=(\delta^{-1}+|\nabla_{y'}\Phi|^2)^{-1}
\bigl( \, (\delta^{-1}-\Delta_{y'})e^{i\Phi}+
i\Delta_{y'}\Phi e^{i\Phi}\, \bigr).
$$
Consequently, if we integrate by parts using this formula, we conclude
that, for a given $N$, we can write $K_\delta$ as a finite combination
of terms of the form
$$\int e^{i\Phi}\, \frac{\delta^{-l}(\partial/\partial y')^{\alpha_1}\Phi\dots
(\partial/\partial
y')^{\alpha_m}\Phi}{(\delta^{-1}+|\nabla_{y'}\Phi|^2)^{j+k+l}}
(\partial/\partial y')^{\gamma} a^2_\delta \, dy,
\tag 3.14$$
where $j+k+l=2N$, $m\le j$, and $|\gamma|\le 2k$.  Using $(3.12)$ and $(3.13)$
we conclude that
$$
(\delta^{-1}+|\nabla_{y'}\Phi|^2)^{-1}\, |(\partial/\partial
y')^\alpha\Phi|\le C_\alpha \delta^{-1/2}
(\delta^{-1}+|\xi'-\eta'|^2)^{-1/2}.$$
Since 
$$(\partial/\partial y')^\gamma a^2_\delta =
O(\delta^{-|\gamma|/2})=O(\delta^{-k}),$$
we conclude that $(3.14)$ is majorized by
$$\delta^{(n-1)/2}\delta^{-j/2-k-l}
(\delta^{-1}+|\xi'-\eta'|^2)^{-j/2-k-l},$$
yielding the other bound for $K_\delta$, which finishes the proof.  \qed
\enddemo

So far we have shown that (3.7) holds when $0<\delta<\delta_1$, with
$\delta_1$ being a uniform small constant.  The argument for the
remaining case where $\delta_1<\delta<10$ is easier.  We note that if
$t-|x|\ge \delta_1>0$ then the above arguments show that if $b\in
C^\infty$
vanishes near $\xi=0$ but equals $1$ outside of a sufficiently small
neighborhood of the origin, then for $\text{Re }z=0$
$$\align
&(S_zg)(t,x)
\\
&\quad=e^{z^2/2}\iint_{1\le|y|\le2}e^{i(x-y)\cdot\xi-i(t-|y|)|\xi|}|\xi|^{-z}
\rho(t^{1-\alpha}\delta^\alpha\xi)
b(y/|y|-\xi/|\xi|)g(y)dyd\xi 
\\
&\quad+O(t^{-N}),
\endalign$$
for any $N$.  If we call $\tilde S_zg$ the first term on the right,
then we need only estimate it.  By Plancherel's theorem
$$\|\tilde S_zg(t,\cd)\|_{L^2}^2\le C\int\Bigl|\, \int_{1\le
|y|\le2}e^{iy\cdot\xi
-i|y|\, |\xi|}b(y/|y|-\xi/|\xi|)g(y) dy\, \Bigr|^2\, d\xi.$$
Since $\text{det }\partial^2\phi/\partial y_j\partial \xi_k\ne 0$ on
the support of the symbol, where $\phi=y\cdot\xi-|y|\, |\xi|$ is the
phase, we can use H\"ormander's $L^2$ theorem for Fourier integral
operators to conclude that the last term is dominated by
$\|g\|_{L^2}^2$.  From this we conclude that (3.7) must hold when
$\delta>\delta_1$,
which finishes our proof.

\head{\bf 4. $L^2$ estimates}\endhead

To finish matters, we still have to prove $(1.13)$.  Since it is easy
to handle small times, we see that it suffices to show that if $w$ solves
the inhomogeneous wave equation $\square w=F$ with zero data, and if
$F(t,x)=0$ when $t-|x|\le 1$, then for $T\ge 10$, say,
$$\|\wght^{-1/2-\epsilon}w\|_{L^2(\{(t,x): \, \, T/2\le t\le T\})}
\le CT^{-\epsilon/4}\|\wght^{1/2+\epsilon}F\|_{L^2}.$$
If we split $w$ up as before, $w=w^0+w^1$, where $\square w^1=F^1$, 
with $F^1(t,x)=F(t,x)$ for $t>T/10$ and $0$ otherwise, then it suffices
to show that for $j=0,1$
$$\|\wght^{-1/2-\epsilon}w^j\|_{L^2(\{(t,x): \, \, T/2\le t\le T\})}
\le CT^{-\epsilon/4}\|\wght^{1/2+\epsilon}F^j\|_{L^2}.
\tag4.1$$

Like before, the estimate for $j=1$ is the easiest.  If we repeat the arguments
which showed how (2.1) implies (2.5), we conclude that the version of 
(4.1) for $j=1$ would be a consequence of the following variant of
(2.1) where $w$ and $F$ are now assumed to be as in Proposition 2.1:
$$\|\wght^{-1/2}w\|_{L^2(|x|<t/2, \, \, T/2\le t\le T)}
\le C(\log T)^{1/2}\|\wght^{1/2}F\|_{L^2}.
\tag4.2$$
However, since the proof of (2.4) also shows that the same estimate holds
when $\text{Re }z=1$, we obtain (4.2) and hence (4.1) when $j=1$.

To handle the case where $j=0$, notice first that
the arguments from the preceding section imply that the remaining 
case of (4.1) would follow from showing that if
$$\supp G\subset \{\, (t,x): \, \, 1\le t\le 2, \, \, \delta_0
\le t-|x|\le 2\delta_0\,  \},$$
and if $\square v=G$ with zero data, then for $T\ge 10$ and $\delta\ge
\delta_0$
$$T^{-1/2-\epsilon/2}\delta^{-1/2+\epsilon/2}
\|v\|_{L^2(\{(t,x): \, \, T/2\le t\le T, \, \, 
\delta\le t-|x|\le 2\delta\, \})}
\le C\delta_0^{1/2}\|G\|_{L^2}.
\tag4.3$$
As in \S 3, the case where $t-|x|>10$ is easy to handle using pointwise
estimates for $E_+(t-s,x-y)$ for such $(t,x)$ if $(s,y)\in \supp G$.  So
in what follows we shall assume that $\delta_0\le \delta\le 10$.

To prove (4.3) for $t-|x|\le 10$, it is convenient to split $v$ into
a low and high frequency part.  To this end, fix $\beta\in C^\infty_0(\Rn)$
satisfying $\beta=1$ near the origin.  If we then let $v=v_0+v_1$, where
$$v_0=\iint e^{i(x-y)\cdot\xi}\beta(\delta\xi)\sin((t-s)|\xi|)G(s,y)ds dy d\xi/|\xi|,$$
then it suffices to show that (4.3) holds when $v$ is replaced by $v_j$, $j=0,1$.
Since 
$$(1-\beta(\delta\xi))/|\xi|=O(\delta),$$ 
the bound for the 
high frequency part follows
from Schwarz's inequality and the variant of (3.5) where $|\xi|^{-z}$, $\text{Re }
z=0$, is replaced by $\delta^{-1}(1-\beta(\delta\xi))/|\xi|$.  Since this inequality
follows from the proof of (3.5), we are left with estimating $v_0$.

For this piece, let us notice that 
$$\int_{|\xi|\le1} e^{i(x-y)\cdot\xi}\beta(\delta\xi)\sin((t-s)|\xi|) d\xi/|\xi|
=O((1+|x-y|)^{-(n-1)/2}).$$
Based on this, we conclude that the variant of (4.3) holds if we replace $v$
by
$$\iint_{|\xi|\le1}e^{i(x-y)\cdot\xi}\beta(\delta\xi)\sin((t-s)|\xi|)G(s,y)
ds dy d\xi/|\xi|.$$
Consequently,
our proof of (4.3) and hence $(1.13)$ would
be complete if we could show that when
$$\tilde v(t,x)
=\iint_{|\xi|\ge1}
e^{i(x-y)\cdot\xi+i(t-s)|\xi|}|\xi|^{-1}\beta(\delta\xi)G(s,y)
ds dy d\xi,$$
we have
$$T^{-1/2}\delta^{-1/2}\|\tilde v\|_{L^2(\{(t,x): \, \, T/2\le t\le T, \, \, 
\delta\le t-|x|\le 2\delta\, \})}
\le C(1+|\log\delta|) \delta_0^{1/2}\|G\|_{L^2}.
\tag4.4$$
Here we are assuming that $G$ is above.  Also, notice that the bounds here
are stronger than those in (4.3).

The first step in proving (4.4) is to notice that the Schwarz inequality and
Lemma 3.2 yield
$$\align
&T^{-1/2}\delta^{-1/2}\|\tilde v\|_{L^2(\{(t,x): \, \, T/2\le t\le T, \, \, 
\delta\le t-|x|\le 2\delta\, \})}
\\
&\qquad\le C\sum_{k=0}^\infty
\Bigl(\int\Bigl|\iint_{2^k\le |\xi|\le 2^{k+1}}
e^{i(x-y)\cdot\xi-is|\xi|}|\xi|^{-1/2}\beta(\delta\xi)G(s,y)d\xi ds dy\Bigr|^2
dx\Bigr)^{1/2}
\endalign$$
Next, if we recall the support properties of $G$ and
use Schwarz's inequality as before we find that the right side is
dominated by $\delta_0^{1/2}$ times
$$\sum_{k=0}^\infty \Bigl(\iint
\Bigr|\iint_{2^k\le |\xi|\le 2^{k+1}}e^{i(x-y)\cdot\xi-i(\tau+|y|)|\xi|}
|\xi|^{-1/2}\beta(\delta\xi)G(\tau+|y|,y)d\xi dy\Bigr|^2dxd\tau\Bigr)^{1/2}.
$$
Notice that the $k$-th summand vanishes 
if $k$ is larger than a fixed multiple of $(1+|\log\delta|)$ 
since $\beta\in C^\infty_0$.  Therefore,
if we now apply the dual version of Lemma 3.2, we obtain (4.4).

This completes the proof of $(1.13)$.  \qed

\bigskip

\subhead
Related Estimates
\endsubhead

The above arguments can also be used to prove weighted $L^2$ estimates
for operators which are similar to the solution operator for the inhomogeneous
wave equation with zero Cauchy data $\square w =F$.  As noted before, this equation
is solved via $w=E_+*F$, where 
$E_+(t,x)=\pi^{(1-n)/2}/2\cdot \chi^{-(n-1)/2}_+(t^2-|x|^2)$ for $t\ge0$ and
$0$ otherwise.  

We could also, as in \cite{24}, consider the related analytic family of operators
$$(T^zF)(t,x)=e^{z^2}\int_0^t \int_{\Rn}
\chi^z_+((t-s)^2-|x-y|^2)F(s,y) \, dy ds,$$
where the convolution is interpreted in the sense of distributions.  If 
$\text{Re } z\ge -(n+1)/2$, recall that $T^z: \, L^2_{\text{comp}}(\st)
\to L^2_{\text{loc}}(\st)$.  As a key step in the proof of his estimates,
Strichartz \cite{24} showed that for the critical values $\text{Re } z=-(n+1)/2$,
$T^z: \, L^2(\st) \to L^2(\st)$.

The above arguments show that a weighted version of this estimate holds under
our support assumptions.  Specifically, if we assume that $F(t,x)=0$ when
$t-|x|\le 0$ or $t<0$ and if $\epsilon>0$, then
$$\multline
\|\, \wght^{-\epsilon+((n+1)/2+\sigma)/2}T^zF\, \|_{L^2(\st)}
\\
\le C\|\, \wght^{\epsilon-((n+1)/2+\sigma)/2}F\|_{L^2(\st)}, 
\endmultline
\tag4.5$$
provided that
$$-(n+1)/2\le \sigma=\text{Re } z\le -(n-1)/2.$$
Georgiev \cite{2} showed how this estimate along with a natural extension
of John's \cite{6} $L^\infty$ estimates can be used to prove non-trivial
weighted estimates off of the line of duality.  Further details will be given
later.

By Stein's analytic interpolation theorem, to prove (4.5), it suffices to
handle the extreme cases where $\text{Re } z=-(n-1)/2$ or $-(n+1)/2$.  The first
case of course follows from the arguments given in this section since,
as we noted before, $T^zF$ behaves essentially like the solution of the inhomogeneous
wave equation when $\text{Re } z=-(n-1)/2$.  Also, since $T^{(n+1)/2+z}$ essentially
agrees with the operator $W^z$ in (2.4), our arguments also yield (4.5)
for the other extreme case where $\text{Re } z=-(n+1)/2$.

\enddocument